\documentclass{elsarticle}
\usepackage{amsfonts}
\usepackage{amssymb}
\usepackage{graphicx}
\usepackage{pstricks}
\usepackage{amsmath}
\usepackage{amsxtra}
\usepackage[english,francais]{babel}

\newtheorem{theorem}{Theorem}[section]
\newtheorem{lemma}[theorem]{Lemma}
\newtheorem{problem}[theorem]{Problem}
\newtheorem{proposition}[theorem]{Proposition}
\newtheorem{corollary}[theorem]{Corollary}
\newtheorem{e-definition}[theorem]{Definition\rm}
\newtheorem{remark}[theorem]{Remark}

\def\og{\leavevmode\raise.3ex\hbox{$\scriptscriptstyle\langle\!\langle$~}}
\def\fg{\leavevmode\raise.3ex\hbox{~$\!\scriptscriptstyle\,\rangle\!\rangle$}}

\begin{document}

\begin{frontmatter}


\title{Toral and spherical Aluthge transforms \\ of $2$-variable weighted shifts}


\selectlanguage{english}
\author[authorlabel1]{Ra\'{u}l E. Curto}
\address[authorlabel1]{Department of Mathematics, The University of Iowa, Iowa City, Iowa
52242}
\ead{raul-curto@uiowa.edu}
\author[authorlabel2]{Jasang Yoon}
\address[authorlabel2]{School of Mathematical and Statistical Sciences, The University of Texas
Rio Grande Valley, Edinburg, Texas 78539}
\ead{jasang.yoon@utrgv.edu}



\begin{abstract}
\selectlanguage{english}
We introduce two natural notions of Aluthge transforms (toral and
spherical) for 2-variable weighted shifts and study their basic properties. \ Next, we study the class of spherically quasinormal $2$-variable weighted
shifts, which are the fixed points for the spherical Aluthge transform. \ Finally,
we briefly discuss the relation between spherically quasinormal and
spherically isometric $2$-variable weighted shifts. \
{\it To cite this article: R. Curto, J. Yoon, C. R. Acad. Sci. Paris, Ser. I (2016).}

\vskip 0.5\baselineskip

\selectlanguage{francais}

\noindent{\bf R\'esum\'e} \vskip 0.5\baselineskip \noindent
Nous introduisons deux notions naturelles des transformations d'Aluthge (toral et sph\'erique) pour les shifts pond\'er\'es \`a deux variables et nous \'etudions leurs propri\'et\'es. \ Ensuite, nous \'etudions la classe de shifts pond\'er\'es \`a deux variables sph\'eriques et quasinormaux, qui sont les points fixes pour la transformation d'Aluthge sph\'erique. \ Enfin, nous discutons bri\`evement de la relation entre les shifts pond\'er\'es \`a deux variables qui sont sph\'eriquement quasinormaux, et ceux qui sont sph\'eriquement isom\`etriques. 
{\it Pour citer cet article~: R. Curto, J. Yoon, C. R. Acad. Sci.
Paris, Ser. I (2016).}

\end{abstract}
\end{frontmatter}

\selectlanguage{francais}

\section*{Version fran\c{c}aise abr\'eg\'ee}

Pour $T\in \mathcal{B}(\mathcal{H})$, la d\'ecomposition polaire de $T$ est $T \equiv U|T|$, o\'{u} $U$ est une isom\'{e}trie partielle et $|T|:=\sqrt{T^{\ast }T}$. \ 
La \textit{transformation d'Aluthge} de $T$ est l'op\'{e}rateur $\widetilde{T}:=|T|^{\frac{1}{2}}U|T|^{\frac{1}{2}}$. \ 
Cette transformation a \'et\'e considere\'e pour la premi\`ere fois dans \cite{Alu} afin d'\'etudier des op\'erateurs $p$-hyponormal et log-hyponormal. Bref, l'id\'ee derri\`ere la transformation d'Aluthge est de convertir un op\'erateur \`a un autre qui partage avec le premier quelques propri\'et\'es spectrales, mais qui reste plus proche d'\^etre un op\'erateur normal. \ Ces derni\`eres ann\'ees, la transformation d'Aluthge a re\c cu de l'attention consid\'erable.

Dans cet article, nous pr\'esentons d'abord la d\'ecomposition polaire des op\'erateurs born\'es et nous \'etudions deux transformations d'Aluthge (toral et sph\'erique) de shifts \`a deux variables $W_{(\alpha ,\beta )}\equiv (T_{1},T_{2})$. \ Puisque a priori il y a plusieurs notions possibles, nous discutons de deux d\'efinitions plausibles et de propri\'et\'es fondamentales dans la deuxi\`eme partie. \ Notre recherche nous permettra de comparer les deux d\'efinitions quant \`a la fa\c con dont elles g\'en\'eralisent la notion d'une variable. Apr\`es avoir discut\'e de quelques propri\'et\'es fondamentales de chaque transformation d'Aluthge, nous proc\'edons \`a \'etudier les deux transformations dans le cas des shifts pond\'er\'es de deux variables. Nous consid\'erons des sujets tels que la pr\'eservation de la hyponormalit\'e conjointe, la continuit\'e de norme, la quasinormalit\'e sph\'erique et l'isom\'etrie sph\'erique.

\selectlanguage{english} 

\section{Introduction}

\label{} For $T\in \mathcal{B}(\mathcal{H})$, the \textit{polar decomposition} of $T$
is $T \equiv U|T|$, where $U$ is a partial isometry and $|T|:=\sqrt{T^{\ast }T}$. \
The \textit{Aluthge transform} of $T$ is the operator $\widetilde{T}:=|T|^{%
\frac{1}{2}}U|T|^{\frac{1}{2}}$. \ This transform was first considered in \cite%
{Alu} in order to study $p$-hyponormal and $\mathrm{\log }$-hyponormal
operators. \ Roughly speaking, the idea behind the Aluthge transform is to
convert an operator into another operator which shares with the first one
some spectral properties, but it is closer to being a normal operator. \ In
recent years, the Aluthge transform has received substantial attention. \
Jung, Ko and Pearcy proved in \cite{JKP} that $T$ has a nontrivial invariant
subspace if and only if $\widetilde{T}$ does. \ (Since every normal operator
has nontrivial invariant subspaces, the Aluthge transform has a natural
connection with the invariant subspace problem.) \ 

For a weighted shift $%
W_{\alpha }\equiv \mathrm{shift}(\alpha _{0},\alpha _{1},\cdots )$, the
Aluthge transform $\widetilde{W}_{\alpha }$ of $W_{\alpha }$ is also a
unilateral weighted shift, given by $\widetilde{W}_{\alpha }\equiv \mathrm{%
shift}(\sqrt{\alpha _{0}\alpha _{1}},\sqrt{\alpha _{1}\alpha _{2}},\cdots )$
(see \cite{LLY}). \ We recall that $W_{\alpha }$ is hyponormal if
and only if $\alpha _{0}\leq \alpha _{1}\leq \cdots $. \ Thus, if $W_{\alpha
}$ is hyponormal, then $\widetilde{W}_{\alpha }$ is also hyponormal. \
However, the converse is not true in general. \ For example, if $W_{\alpha
}\equiv \mathrm{shift}\left( \frac{1}{2},2,\frac{1}{2},2,\cdots \right) $,
then $W_{\alpha }$ is clearly not hyponormal but $\widetilde{W}_{\alpha
}=U_{+}$ is subnormal. \ (Here and in what follows, $U_{+}$ denotes the
(unweighted) unilateral shift.) \ In \cite{LLY}, S.H. Lee, W.Y. Lee, and the
second-named author showed that for $k\geq 2$, the Aluthge
transform, when acting on weighted shifts, need not preserve $k$%
-hyponormality. \ Finally, G. Exner proved in \cite{Exn} that the Aluthge
transform of a subnormal weighted shift need not be subnormal.

In this article, we first introduce the polar decompositions of a commuting pair of bounded
operators and study two Aluthge transforms (toral and spherical) of $2$%
-variable weighted shifts $W_{(\alpha ,\beta )}\equiv (T_{1},T_{2})$. \
Since a priori there are several possible notions, we discuss two plausible
definitions and their basic properties in Section 2. \ Our research will allow us
to compare both definitions in terms of how well they generalize the
$1$-variable notion. \ After discussing some basic properties of each
Aluthge transform, we proceed to study both transforms in the case of $2$%
-variable weighted shifts. \ We consider such topics as preservation of
joint hyponormality, norm continuity, spherical quasinormality, and spherical
isometries. \ 

We first recall that a commuting pair is subnormal if and only if it is $k$-hyponormal for every
$k\geq 1$ (see \cite[Theorem 2.3]{CLY1}). \ For $i=1,2$, we consider the
polar decomposition $T_{i}=U_{i}\left\vert T_{i}\right\vert $ of $T_{i}$. \
Then, for a $2$-variable weighted shift $W_{(\alpha ,\beta )}\equiv
(T_{1},T_{2})$, we define the \textit{toral} Aluthge transform of $%
W_{(\alpha ,\beta )}$ as follows: $\widetilde{W}_{(\alpha ,\beta )}=(%
\widetilde{T}_{1},\widetilde{T}_{2}):=\left( |T_{1}|^{\frac{1}{2}%
}U_{1}|T_{1}|^{\frac{1}{2}},|T_{2}|^{\frac{1}{2}}U_{2}|T_{2}|^{\frac{1}{2}%
}\right) $. \ As we will see in Proposition \ref{Prop1}, the commutativity
of $\widetilde{W}_{(\alpha ,\beta )}$ does not automatically follow from the
commutativity of $W_{(\alpha ,\beta )}$. \ 

There is a second plausible definition of Aluthge transform, which uses a joint polar
decomposition. \ Assume that we have a decomposition of the form $%
(T_{1},T_{2})=\left( U_{1}P,U_{2}P\right) $, where $P:=\sqrt{T_{1}^{\ast
}T_{1}+T_{2}^{\ast }T_{2}}$ and $\textrm{ker } U_1 \cap \textrm{ker } U_2 = \textrm{ker } P$. \ Now let $\widehat{W}_{(\alpha ,\beta )}=(%
\widehat{T}_{1},\widehat{T}_{2}):=\left( \sqrt{P}U_{1}\sqrt{P},\sqrt{P}U_{2}%
\sqrt{P}\right) $. \ We refer to $\widehat{W}_{(\alpha ,\beta )}$ as the
\textit{spherical} Aluthge transform of $W_{(\alpha ,\beta )}$. \ Even though $%
\widehat{T}_{1}=\sqrt{P}U_{1}\sqrt{P}$ is not the Aluthge transform of $%
T_{1} $, we observe that $Q:=\sqrt{U_{1}^{\ast }U_{1}+U_{2}^{\ast }U_{2}}$ is a
(joint) partial isometry; for, $PQ^2P=P^2$, from which it follows that $Q$ is isometric on the range of $P$. \ 

We will prove in
Section 2 that this particular definition of the Aluthge transform
preserves commutativity, and it also behaves well in terms of
hyponormality, for a large class of $2$-variable weighted shifts. \
There is also another useful aspect of the spherical Aluthge transform, which
we now mention. \ If we consider the fixed points of this transform acting
on $2$-variable weighted shifts, we obtain an appropriate generalization of
the concept of quasinormality. \ More precisely, if a $2$-variable weighted
shift $W_{(\alpha ,\beta )}=(T_{1},T_{2})$ satisfies $W_{(\alpha ,\beta )}=(%
\widehat{T}_{1},\widehat{T}_{2})$, then $T_{1}^{\ast }T_{1}+T_{2}^{\ast
}T_{2}$ is, up to scalar multiple, a spherical isometry. \ (We recall that a commuting pair $\mathbf{T} \equiv (T_1,T_2)$ is
called a {\it spherical isometry} if $P^{2}\equiv T_1^*T_1+T_2^*T_2=I$ \cite{EsP}.) \ It follows that
we can then study some properties of the spherical Aluthge transform using well
known results about spherical isometries. \ In this paper, we also focus on the
following basic problem. \smallskip

\begin{problem}
\label{problem 1} (i) For $k\geq 1$, if $W_{(\alpha ,\beta )}$ is $k$%
-hyponormal, does it follow that the toral Aluthge transform $\widetilde{W}%
_{(\alpha ,\beta )}$ $k$-hyponormal? \ What about the case of the spherical
Aluthge transform? \ When does either Aluthge transform preserve hyponormality?; \newline (ii) When do we have the equality 
$\widetilde{W}_{(\alpha ,\beta )}=\widehat{W}_{(\alpha ,\beta )}$?; \newline
(iii) Is the toral Aluthge transform $\mathbf{(}T_{1},T_{2})\rightarrow
\left( \widetilde{T}_{1},\widetilde{T}_{2}\right) $ continuous in the
uniform topology? Similarly, does continuity hold for the spherical Aluthge
transform?
\end{problem}

\section{Main Results}
We first consider (joint) hyponormality for the
toral and spherical Aluthge transforms. \smallskip

\begin{proposition}
\label{Prop1}(i) For $W_{(\alpha ,\beta )}$, we have: (a) $\widetilde{W}%
_{(\alpha ,\beta )}$ is a commuting pair if and only if%
\[
\alpha _{(k_{1},k_{2}+1)}\alpha _{(k_{1}+1,k_{2}+1)}=\alpha
_{(k_{1}+1,k_{2})}\alpha _{(k_{1},k_{2}+2)}\text{ for all }k_{1},k_{2}\geq 0
\]
(with similar conditions holding for the weight sequence $\{\beta _{(k_{1},k_{2})}\}$); and (b)\ $\widehat{W}_{(\alpha ,\beta )}$ is always
a commuting pair. \newline (ii) If $W_{(\alpha ,\beta )}$ is a commuting pair of hyponormal operators, so is $\widehat{W}_{(\alpha
,\beta )}$.
\end{proposition}

\medskip
We next show that there exists a subnormal $W_{(\alpha ,\beta )}$ such that $\widetilde{W}_{(\alpha ,\beta )}$ is not hyponormal; we also prove that there exists a non-hyponormal $W_{(\alpha ,\beta )}$ such that $\widehat{W}_{(\alpha
,\beta )}$ is hyponormal. \ We start with some
definitions. \ The \textit{core} $c(W_{(\alpha ,\beta )})$ of $W_{(\alpha
,\beta )}$ is the restriction of $W_{(\alpha ,\beta )}$ to the invariant
subspace $\mathcal{M}\cap \mathcal{N}$, where $\mathcal{M}\equiv \bigvee
\left\{ e_{(k_{1},k_{2})}\in \ell ^{2}(\mathbb{Z}_{+}^{2}):k_{1}\geq 0\text{
and }k_{2}\geq 1\right\} $ and $\mathcal{N}\equiv \bigvee \left\{
e_{(k_{1},k_{2})}\in \ell ^{2}(\mathbb{Z}_{+}^{2}):k_{1}\geq 1\text{ and }%
k_{2}\geq 0\right\} $. \newline Given a $1$-variable unilateral weighted shift $W_{\omega}$, consider the $2$-variable weighted shift $\Theta
\left( W_{\omega }\right) \equiv W_{(\alpha ,\beta )}$ on $\ell ^{2}(\mathbb{%
Z}_{+}^{2})$ given by the double-indexed weight sequences $\alpha
_{(k_{1},k_{2})}=\beta _{(k_{1},k_{2})}:=\omega _{k_{1}+k_{2}}$ for $%
k_{1},k_{2}\geq 0$. \ The shift $\Theta \left( W_{\omega
}\right) $ can be represented by the weight diagram in Figure \ref{Figure 1}(ii). \smallskip

\begin{lemma}
\label{propscaling2}(\cite{CLY7}) \ Let $\Theta \left( W_{\omega
}\right) $ be given by Figure \ref{Figure 1}(ii), and let $k\geq 1$. \ $W_{\omega }$ is $k$-hyponormal if and only if $\Theta \left( W_{\omega
}\right)$ is (jointly) $k$-hyponormal.
\end{lemma} \smallskip

We next show that hyponormality is not stable under the toral Aluthge transform. \smallskip
\setlength{\unitlength}{1mm} \psset{unit=1mm}
\begin{figure}[th]
\begin{center}
\begin{picture}(162,36)

\psline{->}(10,6)(58,6)
\psline(10,16)(57,16)
\psline(10,26)(57,26)
\psline(10,36)(57,36)
\psline{->}(10,6)(10,38)
\psline(25,6)(25,37)
\psline(40,6)(40,37)
\psline(55,6)(55,37)

\put(2.3,3.2){\footnotesize{$(0,0)$}}
\put(9,-1){$(i)$}
\put(60,-1){$(ii)$}
\put(103,-1){$(iii)$}
\put(21,3){\footnotesize{$(1,0)$}}
\put(36,3){\footnotesize{$(2,0)$}}
\put(51,3){\footnotesize{$(3,0)$}}

\put(14,7){\footnotesize{$\alpha_{(0,0)}$}}
\put(29,7){\footnotesize{$\alpha_{(1,0)}$}}
\put(44,7){\footnotesize{$\alpha_{(2,0)}$}}
\put(56,7){\footnotesize{$\cdots$}}

\put(14,17){\footnotesize{$\alpha_{(0,1)}$}}
\put(29,17){\footnotesize{$\omega_{0}$}}
\put(44,17){\footnotesize{$\omega_{1}$}}
\put(56,17){\footnotesize{$\cdots$}}

\put(14,27){\footnotesize{$\alpha_{(0,2)}$}}
\put(29,27){\footnotesize{$\omega_{1}$}}
\put(44,27){\footnotesize{$\omega_{2}$}}
\put(56,27){\footnotesize{$\cdots$}}

\put(15,36.5){\footnotesize{$\cdots$}}
\put(30,36.5){\footnotesize{$\cdots$}}
\put(45,36.5){\footnotesize{$\cdots$}}

\psline{->}(25,1)(40,1)
\put(30,-2){$\rm{T}_1$}
\psline{->}(7,20)(7,33)
\put(2,27){$\rm{T}_2$}

\put(2.7,15){\footnotesize{$(0,1)$}}

\put(2.7,35){\footnotesize{$(0,3)$}}

\put(10,11){\footnotesize{$\beta_{(0,0)}$}}
\put(10,20){\footnotesize{$\beta_{(1,0)}$}}
\put(10,30){\footnotesize{$\beta_{(2,0)}$}}

\put(25,11){\footnotesize{$\beta_{(0,1)}$}}
\put(25,20){\footnotesize{$\omega_{0}$}}
\put(25,30){\footnotesize{$\omega_{1}$}}

\put(40,11){\footnotesize{$\beta_{(0,2)}$}}
\put(40,20){\footnotesize{$\omega_{1}$}}
\put(40,30){\footnotesize{$\omega_{2}$}}


\psline{->}(62,6)(101,6)
\psline(62,16)(100,16)
\psline(62,26)(100,26)
\psline(62,36)(100,36)
\psline{->}(62,6)(62,38)
\psline(74,6)(74,37)
\psline(86,6)(86,37)
\psline(98,6)(98,37)

\put(60,3.2){\footnotesize{$(0,0)$}}
\put(70,3){\footnotesize{$(1,0)$}}
\put(82,3){\footnotesize{$(2,0)$}}
\put(95,3){\footnotesize{$(3,0)$}}

\put(68,7){\footnotesize{$\omega_{0}$}}
\put(79,7){\footnotesize{$\omega_{1}$}}
\put(91,7){\footnotesize{$\omega_{2}$}}
\put(99,7){\footnotesize{$\cdots$}}

\put(68,17){\footnotesize{$\omega_{1}$}}
\put(79,17){\footnotesize{$\omega_{2}$}}
\put(91,17){\footnotesize{$\omega_{3}$}}
\put(99,17){\footnotesize{$\cdots$}}

\put(68,27){\footnotesize{$\omega_{2}$}}
\put(79,27){\footnotesize{$\omega_{3}$}}
\put(91,27){\footnotesize{$\omega_{4}$}}
\put(99,27){\footnotesize{$\cdots$}}

\put(68,36.5){\footnotesize{$\cdots$}}
\put(79,36.5){\footnotesize{$\cdots$}}
\put(91,36.5){\footnotesize{$\cdots$}}

\psline{->}(75,1)(90,1)
\put(80,-2){$\rm{T}_1$}

\put(62,10){\footnotesize{$\omega_{0}$}}
\put(62,20){\footnotesize{$\omega_{1}$}}
\put(62,30){\footnotesize{$\omega_{2}$}}

\put(74,10){\footnotesize{$\omega_{1}$}}
\put(74,20){\footnotesize{$\omega_{2}$}}
\put(74,30){\footnotesize{$\omega_{3}$}}

\put(86,10){\footnotesize{$\omega_{2}$}}
\put(86,20){\footnotesize{$\omega_{3}$}}
\put(86,30){\footnotesize{$\omega_{4}$}}


\psline{->}(105,6)(153,6)
\psline(105,16)(152,16)
\psline(105,26)(152,26)
\psline(105,36)(152,36)
\psline{->}(105,6)(105,38)
\psline(120,6)(120,37)
\psline(135,6)(135,37)
\psline(150,6)(150,37)

\put(103,3.2){\footnotesize{$(0,0)$}}
\put(112,3){\footnotesize{$(1,0)$}}
\put(131,3){\footnotesize{$(2,0)$}}
\put(146,3){\footnotesize{$(3,0)$}}

\put(109,7){\footnotesize{$\sqrt{\omega_{0}\omega_{1}}$}}
\put(124,7){\footnotesize{$\sqrt{\omega_{1}\omega_{2}}$}}
\put(139,7){\footnotesize{$\sqrt{\omega_{2}\omega_{3}}$}}
\put(151,7){\footnotesize{$\cdots$}}

\put(109,17){\footnotesize{$\sqrt{\omega_{1}\omega_{2}}$}}
\put(124,17){\footnotesize{$\sqrt{\omega_{2}\omega_{3}}$}}
\put(139,17){\footnotesize{$\sqrt{\omega_{3}\omega_{4}}$}}
\put(151,17){\footnotesize{$\cdots$}}

\put(109,27){\footnotesize{$\sqrt{\omega_{2}\omega_{3}}$}}
\put(124,27){\footnotesize{$\sqrt{\omega_{3}\omega_{4}}$}}
\put(139,27){\footnotesize{$\sqrt{\omega_{4}\omega_{5}}$}}
\put(151,27){\footnotesize{$\cdots$}}

\put(110,36.5){\footnotesize{$\cdots$}}
\put(125,36.5){\footnotesize{$\cdots$}}
\put(140,36.5){\footnotesize{$\cdots$}}

\put(127,-2.0){$\rm{T}_1$}
\psline{->}(121,1)(136,1)

\put(105,11){\footnotesize{$\sqrt{\omega_{0}\omega_{1}}$}}
\put(105,21){\footnotesize{$\sqrt{\omega_{1}\omega_{2}}$}}
\put(105,31){\footnotesize{$\sqrt{\omega_{2}\omega_{3}}$}}

\put(120,11){\footnotesize{$\sqrt{\omega_{1}\omega_{2}}$}}
\put(120,21){\footnotesize{$\sqrt{\omega_{2}\omega_{3}}$}}
\put(120,31){\footnotesize{$\sqrt{\omega_{3}\omega_{4}}$}}

\put(135,11){\footnotesize{$\sqrt{\omega_{2}\omega_{3}}$}}
\put(135,21){\footnotesize{$\sqrt{\omega_{3}\omega_{4}}$}}
\put(135,31){\footnotesize{$\sqrt{\omega_{4}\omega_{5}}$}}

\end{picture}
\end{center}
\caption{Weight diagram of the $2$-variable weighted shift with $c(W_{(\protect%
\alpha ,\protect\beta )})=\Theta \left( W_{\protect\omega }\right)$, weight
diagram of a generic $2$-variable weighted shift $\Theta \left( W_{\protect%
\omega }\right) \equiv \mathbf{(}T_{1},T_{2})$, and weight diagram of the
Aluthge transform $\widetilde{\Theta }\left( W_{\protect\omega }\right)
= \widehat{\Theta }\left( W_{\protect\omega }\right)$ of $\Theta
\left( W_{\protect\omega}\right)$, respectively.}
\label{Figure 1}
\end{figure}
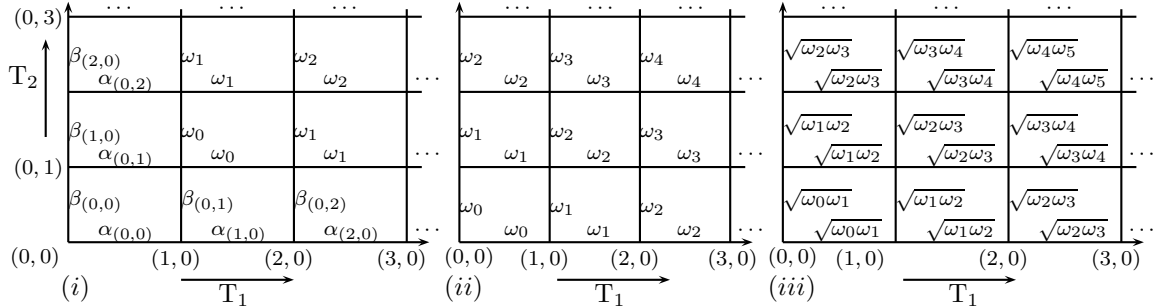
\smallskip

\begin{proposition}
\label{2 atomic}  For $0<x,y<1$, let $W_{(\alpha ,\beta )}$ be the $2$-variable weighted shift in Figure \ref{Figure 1}(i), where $\omega_0=\omega_1=\omega_2=\cdots:=1$, $\alpha_{(0,0)}=\beta_{(0,0)}:=x$, $\alpha_{(0,k_2)}=\beta_{(k_1,0)}:=y \; (k_1,k_2 \ge 1)$ and all remaining weights are equal to $1$. \ Then, we have: \newline
(i) \ $\ W_{(\alpha ,\beta )}$ is subnormal if and only if $x \le s(y):=\sqrt{\frac{1}{2-y^2}}$; \newline
(ii) \ $\ W_{(\alpha ,\beta )}$ is hyponormal if and only if $x\le h(y):=\sqrt{\frac{%
1+y^{2}}{2}}$; \newline
(iii) \ $\ \widetilde{W}_{(\alpha ,\beta )}$ is hyponormal if and only if $x \le CA(y):=\frac{1+y}{2}$; \newline
(iv) \ $\widehat{W}_{(\alpha ,\beta )}$ is hyponormal if and only if $x\le
PA(y):=\frac{2\left( 1+y^{2}-y^{4}\right) }{\left( 1+\sqrt{2}\right) \left(
1+y^{2}\right) \left( \sqrt{1+y^{2}}-y^{2}\right) }$.\newline
Clearly, $s(y) \le h(y) \le PA(y)$ and $CA(y) \le h(y)$ for all $0 < y <1$, while $CA(y) < s(y)$ on $(0,q)$ and $CA(y) > s(y)$ on $(q,1)$, where $q \cong 0.52138$. \ Thus, $W_{(\alpha ,\beta )}$ is hyponormal but $\widetilde{W}%
_{(\alpha ,\beta )}$ is not hyponormal if $0<CA(y)<x\le s(y)$, and $\widehat{W}_{(\alpha ,\beta )}$ is hyponormal but $%
W_{(\alpha ,\beta )}$ is not hyponormal if $0<h(y)<x\le PA(y)$.
\end{proposition}

\smallskip
Next, we consider the invariance of hyponormality under the two Aluthge
transforms. \ We describe a large class of $2$-variable weighted shifts for which both transforms coincide. \ By Figure \ref{Figure 1}(ii) and (iii), and direct calculations, we observe that $\widetilde{\Theta }\left( W_{\omega }\right) =\widehat{\Theta }\left( W_{\omega }\right) $. \ We now have: 
\smallskip

\begin{theorem}
\label{pre-hypo}If $\Theta \left( W_{\omega }\right)$ is hyponormal,
then $\widetilde{\Theta }\left( W_{\omega }\right)$ and $\widehat{\Theta }\left(
W_{\omega }\right)$ are both hyponormal.
\end{theorem}
\smallskip

\begin{remark}
As in the $1$-variable case, we use $W_{\omega }\equiv \mathrm{shift}\left( \frac{1}{2},2,%
\frac{1}{2},2,\cdots \right) $ to build an example of a $2$-variable weighted shift which is not subnormal, but whose Aluthge transforms are both subnormal.
\end{remark}
\smallskip

\setlength{\unitlength}{1mm} \psset{unit=1mm}
\begin{figure}[th]
\begin{center}
\begin{picture}(145,40)

\psline{->}(10,7)(47,7)
\psline(10,17)(45,17)
\psline(10,27)(45,27)

\psline{->}(10,8)(10,35)
\psline(20,7)(20,34)
\psline(30,7)(30,34)
\psline(40,7)(40,34)

\put(8,-1){$(i)$}
\put(59,-1){$(ii)$}
\put(110,-1){$(iii)$}

\put(3,4.5){\footnotesize{$(0,0)$}}
\put(17,3.5){\footnotesize{$(1,0)$}}
\put(27,3.5){\footnotesize{$(2,0)$}}
\put(37,3.5){\footnotesize{$(3,0)$}}

\put(15,7.5){\footnotesize{$a$}}
\put(25,7.5){\footnotesize{$b$}}
\put(35,7.5){\footnotesize{$c$}}
\put(45,7.5){\footnotesize{$d$}}

\put(15,17.5){\footnotesize{$b$}}
\put(25,17.5){\footnotesize{$c$}}
\put(35,17.5){\footnotesize{$d$}}
\put(45,17.5){\footnotesize{$e$}}

\put(15,27.5){\footnotesize{$c$}}
\put(25,27.5){\footnotesize{$d$}}
\put(35,27.5){\footnotesize{$e$}}
\put(45,27.5){\footnotesize{$f$}}

\psline{->}(20,1.5)(30,1.5)
\put(23,-2){${T}_1$}
\psline{->}(2.5,15)(2.5,25)
\put(-2,20){${T}_2$}

\put(3,16){\footnotesize{$(0,1)$}}
\put(3,26){\footnotesize{$(0,2)$}}

\put(10,11){\footnotesize{$y$}}
\put(10,21){\footnotesize{$\frac{by}{a}$}}
\put(10,31){\footnotesize{$\frac{cy}{a}$}}

\put(20,11){\footnotesize{$\frac{by}{a}$}}
\put(20,21){\footnotesize{$\frac{cy}{a}$}}
\put(20,31){\footnotesize{$\frac{dy}{a}$}}

\put(30,11){\footnotesize{$\frac{cy}{a}$}}
\put(30,21){\footnotesize{$\frac{dy}{a}$}}
\put(30,31){\footnotesize{$\frac{ey}{a}$}}

\put(40,11){\footnotesize{$\frac{dy}{a}$}}
\put(40,21){\footnotesize{$\frac{ey}{a}$}}
\put(40,31){\footnotesize{$\frac{fy}{a}$}}


\psline{->}(60,7)(95,7)
\psline(60,17)(94,17)
\psline(60,27)(94,27)

\psline{->}(60,7)(60,35)
\psline(70,7)(70,34)
\psline(80,7)(80,34)
\psline(90,7)(90,34)

\put(55,3.5){\footnotesize{$(0,0)$}}
\put(67,3.5){\footnotesize{$(1,0)$}}
\put(77,3.5){\footnotesize{$(2,0)$}}
\put(87,3.5){\footnotesize{$(3,0)$}}

\put(65,7.5){\footnotesize{$a$}}
\put(75,7.5){\footnotesize{$b$}}
\put(85,7.5){\footnotesize{$c$}}
\put(95,7.5){\footnotesize{$d$}}

\put(65,17.5){\footnotesize{$b$}}
\put(75,17.5){\footnotesize{$c$}}
\put(85,17.5){\footnotesize{$d$}}
\put(95,17.5){\footnotesize{$e$}}

\put(65,27.5){\footnotesize{$c$}}
\put(75,27.5){\footnotesize{$d$}}
\put(85,27.5){\footnotesize{$e$}}
\put(95,27.5){\footnotesize{$f$}}

\psline{->}(70,1.5)(80,1.5)
\put(73,-2){${T}_1$}

\put(49,22){$=$}

\put(60,11){\footnotesize{$a$}}
\put(60,21){\footnotesize{$b$}}
\put(60,31){\footnotesize{$c$}}

\put(70,11){\footnotesize{$b$}}
\put(70,21){\footnotesize{$c$}}
\put(70,31){\footnotesize{$d$}}

\put(80,11){\footnotesize{$c$}}
\put(80,21){\footnotesize{$d$}}
\put(80,31){\footnotesize{$e$}}

\put(90,11){\footnotesize{$d$}}
\put(90,21){\footnotesize{$e$}}
\put(90,31){\footnotesize{$f$}}


\put(100,22){$\circ$}
\put(96,35){Schur}
\put(95,32){product}
\psline{->}(100.8,31)(100.8,25)

\psline{->}(110,7)(146,7)
\psline(110,17)(146,17)
\psline(110,27)(146,27)

\psline{->}(110,7)(110,35)
\psline(120,7)(120,34)
\psline(130,7)(130,34)
\psline(140,7)(140,34)

\put(105,3.5){\footnotesize{$(0,0)$}}
\put(117,3.5){\footnotesize{$(1,0)$}}
\put(127,3.5){\footnotesize{$(2,0)$}}
\put(137,3.5){\footnotesize{$(3,0)$}}

\put(115,7.5){\footnotesize{$1$}}
\put(125,7.5){\footnotesize{$1$}}
\put(135,7.5){\footnotesize{$1$}}
\put(145,7.5){\footnotesize{$1$}}

\put(115,17.5){\footnotesize{$1$}}
\put(125,17.5){\footnotesize{$1$}}
\put(135,17.5){\footnotesize{$1$}}
\put(145,17.5){\footnotesize{$1$}}

\put(115,27.5){\footnotesize{$1$}}
\put(125,27.5){\footnotesize{$1$}}
\put(135,27.5){\footnotesize{$1$}}
\put(145,27.5){\footnotesize{$1$}}

\psline{->}(120,1.5)(130,1.5)
\put(123,-2){${T}_1$}

\put(110,11){\footnotesize{$\frac{y}{a}$}}
\put(110,21){\footnotesize{$\frac{y}{a}$}}
\put(110,31){\footnotesize{$\frac{y}{a}$}}

\put(120,11){\footnotesize{$\frac{y}{a}$}}
\put(120,21){\footnotesize{$\frac{y}{a}$}}
\put(120,31){\footnotesize{$\frac{y}{a}$}}

\put(130,11){\footnotesize{$\frac{y}{a}$}}
\put(130,21){\footnotesize{$\frac{y}{a}$}}
\put(130,31){\footnotesize{$\frac{y}{a}$}}

\put(140,11){\footnotesize{$\frac{y}{a}$}}
\put(140,21){\footnotesize{$\frac{y}{a}$}}
\put(140,31){\footnotesize{$\frac{y}{a}$}}
\end{picture}
\end{center}
\caption{Weight diagram of the 2-variable weighted shift in Theorem \protect
\ref{Thm1}.}
\label{Figure 2}
\end{figure}

\begin{theorem}
\label{Thm1} For $W_{(\alpha ,\beta )}$, $\widehat{W}_{(\alpha ,\beta )}=\widetilde{W}_{(\alpha ,\beta
)} $ if and only if $W_{(\alpha ,\beta )}$ is given by Figure \ref{Figure 2}%
(i). \ Furthermore, if $\mathrm{shift}\left( a,b,c,\cdots \right) $ is
subnormal, then $W_{(\alpha ,\beta )}$ is subnormal with Berger measure $\mu
=\upsilon \ast \delta _{\left( 1,\sqrt{\frac{y}{a}}\right) }$, where $%
\upsilon $ is the diagonal Berger measure of the $2$-variable weighted shift
given by Figure \ref{Figure 2}(ii) (see \cite{CLY7}) and $\ast $ is the
convolution product (see \cite{KiYo2}).
\end{theorem}
\smallskip

We now turn our attention to the
continuity properties of the toral (resp. spherical) Aluthge transform of a
commuting pair. \ Since the continuity of the toral Aluthge transform is straightforward, we focus on the spherical case. \ The following result is well known:\ for a single operator
$T\in \mathcal{B}(\mathcal{H})$, the Aluthge transform map $T\rightarrow
\widetilde{T}$ is $\left( \left\Vert \cdot \right\Vert ,\left\Vert \cdot
\right\Vert \right) -$ continuous on $\mathcal{B}(\mathcal{H})$ (\cite{DySc}%
). \ We extend this to the multivariable case. \smallskip

\begin{lemma}
\label{Re 4}(cf. \cite[Lemma 2.1]{DySc}) \ Let $\mathbf{T}\equiv (T_{1},T_{2})$ be a pair of commuting operators, and let $(T_{1},T_{2})\equiv \left( U_{1}P,U_{2}P\right) $ be its joint polar decomposition; recall that $P=\sqrt{T_{1}^{\ast }T_{1}+T_{2}^{\ast }T_{2}}$. \ For $n\in \mathbb{N}$ and $t>0$, let $f_{n}\left(
t\right) :=\sqrt{\max \left( \frac{1}{n},t\right) }$ and let $A_n:=f_n(\mathbf{T})$. \ Then we have: \newline
(i) \ $\left\Vert A_{n} \right\Vert \leq \max \left( n^{-%
\frac{1}{2}},\left\Vert P\right\Vert ^{\frac{1}{2}}\right) $; \ (ii) \ $\left\Vert P A_n^{-1}\right\Vert
\leq \left\Vert P\right\Vert ^{\frac{1}{2}}$; \ (iii) \ $\left\Vert A_n -P^{\frac{1}{2}}\right\Vert \leq
n^{-\frac{1}{2}}$; \newline 
(iv) \ $\left\Vert PA_n^{-1}-P^{\frac{1}{2}%
}\right\Vert \leq \frac{1}{4} n^{-\frac{1}{2}}$; \ (v) \ for $i=1,2$, $\left\Vert A_n T_{i}A_n^{-1}-P^{\frac{1}{2}}U_{i} P^{\frac{1}{2}}\right\Vert \leq
\frac{5}{4} n^{-\frac{1}{2}}\left\Vert T_{i}\right\Vert ^{\frac{1}{2}}$.
\end{lemma}
\smallskip

\medskip
From Lemma \ref{Re 4} we obtain:

\smallskip
\begin{lemma}
\label{Re 5}(cf. \cite[Lemma 2.2]{DySc}) \ Given $R\geq 1$ and $\epsilon >0$, there are real polynomials $p$ and $q$ such that for every commuting pair $\mathbf{T}%
\equiv (T_{1},T_{2})$ with $\left\Vert T_{i}\right\Vert
\leq R$ $\left( i=1,2\right) $, we have 
\begin{equation*}
\left\Vert P^{\frac{1}{2}}U_{i}P^{\frac{1}{2}}-p\left( T_{1}^{\ast
}T_{1}+T_{2}^{\ast }T_{2}\right) T_{i}q\left( T_{1}^{\ast }T_{1}+T_{2}^{\ast
}T_{2}\right) \right\Vert <\epsilon \text{.}
\end{equation*}
\end{lemma}
\smallskip

By Lemmas \ref{Re 4} and \ref{Re 5}, we have:

\begin{theorem}
\label{ContinuityP}The spherical Aluthge transform $\mathbf{(}%
T_{1},T_{2})\rightarrow \widehat{\mathbf{(}T_{1},T_{2})}$ is $\left(
\left\Vert \cdot \right\Vert ,\left\Vert \cdot \right\Vert \right) -$
continuous on $\mathcal{B}(\mathcal{H})$.
\end{theorem}

\smallskip

We now study the class of spherically quasinormal
(resp. spherically isometric) commuting pairs of operators (\cite{Ath1}, \cite{AtPo}, \cite{AtLu}, \cite{EsP}, \cite{Gle}, \cite{Gle2}). \ In the literature, spherical quasinormality of a commuting $n$-tuple $\mathbf{T} \equiv (T_1,\cdots,T_n)$ is associated with the commutativity of each $T_i$ with $P^2$. \ It is not hard to prove that, for $2$-variable weighted shifts, this is equivalent to requiring that $W_{(\alpha ,\beta )}\equiv (T_1,T_2)$ be a fixed point of the spherical Aluthge transform, that is, $\widehat{W_{(\alpha ,\beta )}}=W_{(\alpha ,\beta )}$. \ A straightforward calculation shows that this is equivalent to requiring that each $U_i$ commutes with $P$. \ In particular, $(U_1,U_2)$ is commuting whenever $(T_1,T_2)$ is commuting. \ Also, recall from Section 1 that a commuting pair $\mathbf{T}$ is a spherical isometry if $P^{2}=I$. \ Thus, in the case of spherically quasinormal $2$-variable weighted shifts, we always have $U_1^*U_1+U_2^*U_2=I$. \ In the following result, the key new ingredient is the equivalence of (i) and (ii). \smallskip

\begin{theorem}
\label{Quasinormal3} For $W_{(\alpha ,\beta )}\equiv (T_{1},T_{2})$, the following statements are equivalent:\newline
(i) $W_{(\alpha ,\beta )}$ is spherically quasinormal; \ (ii) for all $%
(k_{1},k_{2})\in \mathbb{Z}_{+}^{2}$, $\alpha _{(k_{1},k_{2})}^{2}+\beta
_{(k_{1},k_{2})}^{2}=C>0$; \ (iii) $T_{1}^{\ast }T_{1}+T_{2}^{\ast
}T_{2}=C\cdot I$.
\end{theorem}
\smallskip

\textbf{Sketch of Proof.} Briefly stated, our strategy is as follows: By
the continuous functional calculus, we can assume that $T_{1}$ and $T_{2}$
commute with $P$. \ It follows that for all $(k_{1},k_{2})\in
\mathbb{Z}_{+}^{2}$, $\alpha _{(k_{1},k_{2})}^{2}+\beta _{(k_{1},k_{2})}^{2}$
is constant. \ Next, we compute $T_{1}^{\ast }T_{1}+T_{2}^{\ast }T_{2}$. \qed

\smallskip
By the proof of Theorem \ref{Quasinormal3}, we remark that once the zero-th
row of $T_{1}$, call it $W_{0}$, is given, then the entire $2$-variable
weighted shift is fully determined.
\smallskip

We briefly pause to recall the construction of Stampfli's shift $W_{(\sqrt{a},\sqrt{b},\sqrt{c}%
)^{\wedge }}\equiv \mathrm{shift}\left( (\sqrt{a},\sqrt{b},\sqrt{c})^{\wedge
}\right)$, where $0<\sqrt{a}<\sqrt{b}<\sqrt{c}$. \ From \cite{CuFi1}, $W_{(\sqrt{a},\sqrt{b},\sqrt{c}%
)^{\wedge }}$ is subnormal, with $2$-atomic Berger measure $\xi =\rho
_{0}\delta _{s_{0}}+\rho _{1}\delta _{s_{1}}$, where $\varphi _{0}:=-\frac{%
ab(c-b)}{b-a}$, $\varphi _{1}:=\frac{b(c-a)}{b-a}$, $s_{0}:=\frac{\varphi
_{1}-\sqrt{\varphi _{1}^{2}+4\varphi _{0}}}{2}$, $s_{1}:=\frac{\varphi _{1}+%
\sqrt{\varphi _{1}^{2}+4\varphi _{0}}}{2}$, $\rho _{0}:=\frac{s_{1}-a}{%
s_{1}-s_{0}}$, and $\rho _{1}:=\frac{a-s_{0}}{s_{1}-s_{0}}$. \ We are now ready for 

\medskip
\begin{theorem}
\label{Quasinormal2}Consider a spherically quasinormal $W_{(\alpha ,\beta
)}\equiv (T_{1},T_{2})$ given by Figure \ref{Figure 1}%
(i), where $W_0:=\mathrm{shift}(\alpha _{\left( 0,0\right) },\alpha _{\left(
1,0\right) },\cdots )=W_{(\sqrt{a},\sqrt{b},\sqrt{c})^{\wedge }}$.
For given $\mathbf{k}\equiv (k_{1},k_{2})\in \mathbb{Z}_{+}^{2}$, we let $\alpha
_{(k_{1},k_{2})}^{2}+\beta _{(k_{1},k_{2})}^{2}=\varphi _{1}$, where $\varphi_{1}$ is given as above. \ Then,
$W_{(\alpha ,\beta )}$ is subnormal with Berger measure $\left(
\frac{s_{1}-a}{s_{1}-s_{0}}\right) \delta _{\left( s_{0},s_{1}\right)
}+\left( \frac{a-s_{0}}{s_{1}-s_{0}}\right) \delta _{\left(
s_{1},s_{0}\right) }$.
\end{theorem}
\smallskip

The subnormality of $W_{(\alpha ,\beta )}$ in Theorem \ref{Quasinormal2} is a special case of the following result.
\smallskip
\begin{lemma}
\label{Sp-sub}(\cite{EsP}) \ Any spherical isometry is subnormal.
\end{lemma}
\smallskip
By Theorem \ref{Quasinormal3} and Lemma \ref{Sp-sub}, we obtain: 
\smallskip

\begin{corollary}
\label{Qua-sub}Any spherically quasinormal $2$-variable weighted shift is
subnormal.
\end{corollary}
\smallskip

\begin{remark}
(i) \ A. Athavale and S. Poddar have recently proved that a commuting spherically quasinormal pair is always subnormal \cite[Proposition2.1]{AtPo}; this provides a different proof of Corollary \ref{Qua-sub}. \newline 
(ii) \ In a different direction, let $Q_{\mathbf{T}}(X):=T_1^*XT_1+T_2^*XT_2$. \ By induction, it is easy to prove that if $\mathbf{T}$ is spherically quasinormal, then $Q_{\mathbf{T}}^n(I)=(Q_{\mathbf{T}}(I))^n \; (n \ge 0)$; by \cite[Remark 4.6]{ChSh}, $\mathbf{T}$ is subnormal. 
\end{remark}

\medskip

\textit{Acknowledgments.} \ The authors are deeply indebted to the referee for several helpful comments that improved the presentation. \ Some of the calculations in this paper were obtained using the software tool \textit{Mathematica} \cite{Wol}. \ The first author of this paper was partially supported by NSF Grant DMS-1302666.




\end{document}